\documentclass{article}[12pt]
\usepackage{amsthm,amsfonts,amsmath,amssymb}
\newcommand{\eq}{\begin{equation}}
\newcommand{\en}{\end{equation}}

\newcommand{\eqstar}{\begin{eqnarray*}}
\newcommand{\enstar}{\end{eqnarray*}}

\newcommand{\diff}{{\rm d}}
\newcommand{\prob}{\Bbb P}

\newcommand{\Nat}{\Bbb N}

\newtheorem{theorem}{\large Theorem}
\newtheorem{proposition}[theorem] {\large Proposition}

\newtheorem{lemma}[theorem]{\large Lemma}

\newtheorem{remark}[theorem]{\large Remark}

\begin{document}
\title{Corners and Records of the Poisson Process in Quadrant}
\author{Alexander V. Gnedin\thanks{Postal address:
 Department of Mathematics, Utrecht University,
 Postbus 80010, 3508 TA Utrecht, The Netherlands. email: gnedin@math.uu.nl}}
\date{}
\maketitle

\begin{abstract}
\noindent
The scale-invariant spacings lemma due to  Arratia, Barbour and Tavar{\'e}
establishes the distributional identity of 
a self-similar Poisson process and the set of 
spacings between the points of this process.
In this note we connect this result with properties of a certain set of extreme points
of the unit Poisson process in the
positive quadrant.
\end{abstract}

\vskip0.2cm
\noindent
{\large Keywords}: $k$-records, $k$-corners, self-similar Poisson process, Ignatov's theorem\\ 
\vskip0.2cm
\noindent
\large{2000 Mathematics Subject Classification: Primary 60G70, Secondary 60G55}

\section{Introduction}

For fixed $k>0$ let ${\cal T}^{(k)}$ be the {\it self-similar} (or {\it scale-invariant}) Poisson point process on ${\mathbb R}_+$, 
with intensity function $k/t$.
Let ${\cal S}^{(k)}$ be the point process of {\it spacings} in ${\cal T}^{(k)}$, meaning that the generic point of 
 ${\cal S}^{(k)}$
is a difference $t-s$, where  $t>s$ are some consequitive points of ${\cal T}^{(k)}$
(so $[s,t]\cap{\cal T}^{(k)}=\{s,t\}$). The {\it scale-invariant spacings lemma} due to Arratia, Barbour and Tavar{\'e} 
\cite[Lemma 7.1]{ABT} asserts that

\eq\label{ABT}
{\rm the~ ABT~ lemma:}~~~~~~~~  {\cal S}^{(k)}=_d {\cal T}^{(k)}.
\en 
In this note we re-derive this
remarkable result from the prospective of the theory of records,
and connect it with the circle of ideas around Ignatov's theorem
 \cite{Bunge, BC, BC1, GRogers, Resnick}. 
We choose  the framework of the Poisson point process
in the positive quadrant since this setting is very geometric and allows us to exploit various symmetry properties of the Lebesgue measure.
The connection between $k$-records and $k$-corners in Proposition \ref{thin} 
and the intensity formula (\ref{multdens})
are new.

\par See \cite{ ABT, Arr,  selfsim, GPEw} for other occurrences  of the self-similar Poisson process in combinatorial probability.

\section{Corners and records}

Let $\cal P$ be the  Poisson point process in ${\mathbb R}_+^2$ with unit intensity. 
All point processes considered here have no multiple points, 
a feature which enables us to treat 
these  processes 
as random sets  rather than  counting measures.
We shall interpret an atom $a=(t,x)\in\cal P$ as the value $x$ observed at time $t$.
With probability one no two atoms of $\cal P$ lie on the same vertical or horizontal line,
hence there is a one-to-one correspondence between observation times and observed values.
The coordinate projections will be denoted $\tau(a)=t$ and $\xi(a)=x$.
The process $\cal P$ is locally finite, however this does not apply to its projections:
for every interval $]s,t[\,$  ($0<s<t$) there are infinitely many atoms with $\tau(a)\in \, ]s,t[\,$, and for every 
interval $]x,y[\,$ ($0<x<y$) there are infinitely many atoms with $\xi(a)\in \, ]x,y[\,$.

\par For $k$ a positive integer, a point $a\in {\mathbb R}_+^2$ is said to be a
{\it $k$-corner} of $\cal P$ if 
\begin{itemize}
\item[(I)] either $a\in {\cal P}$ and there are $k-1$ points $b\in {\cal P}$ strictly south-west from $a$,
\item[(II)] or $a\notin {\cal P}$, there are $k-2$ points $b\in {\cal P}$ strictly south-west from $a$,
a point $c\in {\cal P}$ strictly west from $a$ and a point $d\in {\cal P}$ strictly south from $a$.
\end{itemize}
\noindent
To interpret the definition geometrically, suppose a light source allocated at point $a\in{\mathbb R}_+^2$ illuminates the area 
south-west from $a$ including the edges. Generate a rectangular greed, dense
in the quadrant, by drawing all vertical and horizontal lines through
atoms of $\cal P$. The $k$-corners are  the points $a$ of the grid which illuminate exactly $k$ atoms of $\cal P$.

\par We denote ${\cal C}^{(k)}$ the set of $k$-corners and denote ${\cal R}^{(k)}$ its subset defined by the condition (I) alone.
Obviously, ${\cal C}^{(1)}={\cal R}^{(1)}$, but for $k>1$ the inclusion ${\cal R}^{(1)}\subset{\cal C}^{(1)}$ is strict
almost surely.
Following \cite{GRogers} we  call the points $a\in {\cal R}^{(k)}$ {\it $k$-records}.
For $a\in{\cal P}$ the {\it initial rank} of $a$ is one bigger the number of atoms $b\in {\cal P}$ strictly south-west from $a$,
hence a $k$-record is an observation of initial rank $k$.

\par Notably, $\tau({\cal C}^{(k)})=_d \xi({\cal C}^{(k)})$ and $\tau({\cal R}^{(k)})=_d \xi({\cal R}^{(k)})$.
This is seen from the fact that the reflection  $(t,x)\mapsto (x,t)$ about the bisectrix preserves both the coordinatewise partial order 
and the Lebesgue measure, hence preserves the law of $\cal P$.

\par Let $M_t^{(k)}
$ be the $k$-th smallest value observed before $t$, which is the $k$-th minimal point of the Poisson process  
$\{\xi(a):\,a\in {\cal P}\cap ([0,t]\times[0,\infty[)\}$.
It is easily seen that $(M_t^{(k)}, t> 0)$ is a nonincreasing piecewise-constant c{\`a}dl{\`a}g process, whose flats start at the $k$-corners
of $\cal P$. Indeed, $M_t^{(k)}>M_{t-}^{(k)}$ means that there is an atom $(t,x)\in {\cal P}$ 
with $t\in \,\,\tau( \cup_{i\leq k}{\cal R}^{(i)})$; then $x=M_t^{(k)}$ if $x$ is a $k$-record, and 
$M_t^{(k)}=M_{t-}^{(k-1)}$ if the initial rank of $(t,x)$ is less than $k$.
Furthermore, if $(t,x)\in {\cal P}$ is a $j$-record for $j< k$ then $M_t^{(k)}> x$ and $M_s^{(k)}=x$ for $s$ the time of the $(k-j)$th observation
 in $]t,\infty[\,\times [0,x]$.
It follows that $\xi(\cup_{j\leq k}{\cal R}^{(j)})=\xi({\cal C}^{(k)})$, and by symmetry also that
$\tau(\cup_{j\leq k}{\cal R}^{(j)})=\tau({\cal C}^{(k)})$.
However, despite the coincidence of projections, the point processes ${\cal C}^{(k)}$ and $\cup_{j\leq k}{\cal R}^{(k)}$ are very different.

\begin{remark}{\rm 
 The term
 `$k$-record' in the existing literature is ambiguous.
By some authors (see e.g. \cite{Nevzorov}) a $k$-record  is 
a new value of the $k$th minimum caused by
an observation of the initial rank at most $k$,
 and this corresponds to the historically 
first usage of the term in \cite{DK}.
By other authors (especially in the work on Ignatov's theorem, see \cite{Bunge} for a survey) a $k$-record is an observation of the initial 
rank exactly $k$.
According to  
\cite{Nagaraja}, these are $k$-records of types 2 and 1, respectively.
Looking in the earlier work on the order properties of multivariate samples  \cite{Barn}, the 
 $k$-records (in the sense of \cite{Bunge}, or of type 1 in \cite{Nagaraja})
correspond to
 `the $k$th layer 3rd quadrant admissible points'.
Thus our $k$-records are as in \cite{Bunge, Resnick} (hence type 1 in \cite{Nagaraja}), while our $k$-corners are 
the `$k$-records' in the sense of \cite{Nevzorov} (hence type 2 in \cite{Nagaraja}).
}
\end{remark}

\section{Projections and intensity}
 The process $(M_t^{(k)}, t>0)$ is Markovian, with a familiar kind of dynamics \cite{BertoinBook, Entr, BC, BC1, chain}. Given $M_t=x$, the residual life-time in $x$ 
is $E/x$ and the new state when the transition occurs is $Bx$, where the random variables $E$ and $B$ are independent, $E$  is  exponential(1), and
$B$ is beta$(k,1)$ with density
\eq\label{beta}
\prob(B\in \diff z)=kz^{k-1}\diff z\,~~~z\in [0,1].
\en 
The marginal distributions have gamma densities
\eq\label{gamma}
\prob(M_t^{(k)}\in \diff x)=    { e^{-tx}(tx)^{k-1}t\diff x\over \Gamma(k)}\,,~~~~~~~x>0.
\en
A self-similarity property 
\eq\label{selfsim}
(cM_t^{(k)},t>0)=_d (M_{t/c}^{(k)},t>0),~~~~~~c>0
\en
follows from the invariance of the Lebesgue measure under the hyperbolic shifts $(t,x)\mapsto(t/c,cx)$.
The process `enters from the infinity', i.e. has the asymptotic initial value
$M_{0+}^{(k)}=\infty$,  and has the asymptotic
 terminal value 
$M_{\infty-}^{(k)}=0$.

\par In the following result, the assertion about $\tau({\cal R}^{(1)})$ is an instance \cite[Proposition 4.9]{Resnick}, while
the last claim  
is a specialisation of Ignatov's theorem 
in the form of 
\cite[Corollary 5.1]{GRogers}.

\begin{proposition}\label{propo} The point processes $\tau({\cal R}^{(k)})$ for $k=1,2,\ldots$ are iid Poisson, each  with intensity $1/t$.
The process $\tau({\cal C}^{(k)})$ is Poisson with  intensity $k/t$.
The analogous facts are true for $\xi({\cal R}^{(k)})$ and  $\xi({\cal C}^{(k)})$.
\end{proposition}
\proof
Fix $t$ and let $a_1,a_2,\ldots$ be the points of ${\cal P}\cup([0,t]\times [0,\infty[)$ labelled by increase of their $x$-values
 $M^{(1)}_t, M^{(2)}_t, \ldots$.
Because the initial rank of $a_k$ is equal to
$\#\{i:i\leq k,\tau(a_i)\leq\tau(a_k)\},$
the processes $\tau({\cal R}^{(k)})\cap[0,t]$ for $k=1,2,\ldots$ are completely determined by the time projections $(\tau(a_k),k\in \Nat)$, hence they are jointly
independent of the $x$-projections $(M^{(k)}_t, k\in \Nat)$. On the other hand, the initial ranks of observations
after $t$ depend on ${\cal P}\cap([0,t]\times [0,\infty[)$ only through $(M^{(k)}_t, k\in \Nat)$, 
from which follows that the {\it multivariate} point process $(\tau({\cal R}^{(k)}), k\in \Nat)$ on ${\mathbb R}_+$ has independent increments,
meaning that its restrictions to disjoint intervals are independent
(a property also called {\it complete independence} in \cite{Dayley}).
The intensity of each $\tau({\cal R}^{(k)})$ is readily identified as $1/t$ since an observation of initial rank $k$ occurs 
in $[t-\diff t,t]$ precisely when $a_k$ arrives on this interval, and since the law of $\tau(a_k)$ is uniform$[0,t]$.
It follows that each $\tau({\cal R}^{(k)})$ is a self-similar Poisson process with intensity $1/t$.
The multivariate process $(\tau({\cal R})^{(k)}, k\in \Nat)$ is {\it simple},
 hence by a standard result 
from the theory of point processes \cite[p. 205]{Dayley} the component processes $\tau({\cal R}^{(k)})$'s are jointly independent.

\par By symmetry about the bisectrix the above is extended  to  $k$-record values.
Superposing $k$ iid Poisson processes
yields the result about $\tau({\cal C}^{(k)})=\cup_{j\leq k} \tau({\cal R}^{(j)})$,
and finally this is extended to $\xi({\cal C}^{(k)})$ by symmetry.
\endproof

\begin{remark}{\rm
Ignatov's theorem in its classical form asserts that the point processes of $k$-record values, derived from an iid sequence 
(with some distribution function $F$)
 are iid. By application of the probability integral transform, the case of 
arbitrary continuous $F$ is reducible to the instance of $F$ being uniform$[0,1]$.
In its turn, the uniform case is readily covered by Proposition \ref{propo}, because the values of 
the observations in the strip ${\cal P}\cap([0,\infty[\,\times [0,x])$, arranged in their time-order, are iid uniform$[0,1]$.
 For continuous $F$,
this argument for Ignatov's theorem 
 seems to be the shortest known.
Note that the symmetry between record values and record times is lost if the Lebesgue measure ${\rm d}t{\rm d}x$ in the quadrant is replaced
by any other ${\rm d}t\cdot\nu({\rm d}x)$ in ${\mathbb R}_+\times{\mathbb R}$ with nonatomic, sigma-finite $\nu$ such that $\nu[-\infty,x]<\infty$
for $x\in {\mathbb R}$. 
}
\end{remark}

\begin{lemma}\label{mylemma} 
Conditionally given $(M_s^{(k)},s\leq t)$, 
 the vector $(M^{(1)}_t,\ldots,M^{(k-1)}_t)$ is independent of the observation times
$(\tau({\cal R}^{(j)})\cap[0,t], j\in\Nat)$ and has 
the same law as the vector of
$k-1$ minimal order statistics sampled from the uniform distribution on $[0,M_t^{(k)}]$.
\end{lemma}
\proof  Condition on the location of the last $k$-corner before $t$, say $(u,x)$. We have   ${\cal P}\cap (\,]u,t[\,\times [0,x])=\varnothing$.
The rectangular grid spanned on $k$ points involved in the definition of $(u,x)$ is distributed like
the product grid generated by $k-1$ order statistics from uniform$[0,u]$ and  independent $k-1$ order statistics from uniform$[0,x]$.
The grid and the processes ${\cal P}\cap (]t,\infty]\times [0,x])$, ${\cal P}\cap([0,t]\times\,]x,\infty])$ are jointly independent,
which readily yields the result, because 
$(M_s^{(k)},s\leq t)$ is determined by  ${\cal P}\cap([0,u]\times\,]x,\infty])$, and 
$(\tau({\cal R}^{(j)})\cap[0,t], j\in\Nat)$ is determined by ${\cal P}\cap([0,u]\times\,]x,\infty])$ 
and the $\tau$-projection of the grid.
\endproof
\noindent

\begin{remark}\label{remaf}
{\rm 
The $k$th sample from the uniform distribution hits each of the $k$ spacings generated by $k-1$ order statistics 
with probability $1/k$.
Thus 
Lemma \ref{mylemma} implies that each
$\tau({\cal R}^i)$ is a pointwise Bernoulli thinning with probability $1/k$ of the process $\tau({\cal C}^{(k)})$, for any $k\geq i$.
Once the independence of increments of the superposition process $\tau({\cal C}^{(k)})=\cup_{i\leq k}\tau({\cal R}^{(i)})$ is acquired,
these facts can be used
 to avoid the most subtle part of our proof of Proposition \ref{propo}:
the reference to the general result that the independence of increments of a multivariate process and simplicity imply independence
of the marginal point processes.
}
\end{remark}

We compute next the density $p_m(a_1,\ldots,a_m)$ of the event that there are $m$ $k$-corners at locations $a_1=(t_1,x_1),\ldots,a_m=(t_m,x_m)$, 
with $t_1<\ldots<t_m,~x_1>\ldots>x_m$, and no further $k$-corners occur between $t_1$ and $t_m$. 
This event occurs when $M_{t_1-}=x$ for some $x>x_1$ and the process at time $t_1$ decrements to $x_1$, then 
spends the time $t_2-t_1$ at $x_1$, then decrements to $x_2$ and so on, hence the infinitesimal probability of the 
event in focus is
\begin{eqnarray*}
 {1\over\Gamma(k)}\int_{x_1}^\infty (t_1x)^{k-1}e^{-t_1x}t_1\diff x \,(x\diff t_1)\, k\left({x_1\over x}\right)^{k-1}\,
{\diff x_1\over x}\\
x_1e^{-(t_2-t_1)x_1}\diff t_2\,
 k\left({x_2\over x_1}\right)^{k-1}\,{\diff x_2\over x_1}\cdots\\
x_{m-1}e^{-(t_m-t_{m-1})x_{m-1}}\diff t_m\,
 k\left({x_m\over x_{m-1}}\right)^{k-1}\,{\diff x_m\over x_{m-1}},
\end{eqnarray*}
 which after massive cancellation results in
\eq\label{multdens}
p_m(a_1,\ldots,a_m)={k^m\over \Gamma(k)} (t_1 x_m)^{k-1}\exp\big[{-x_1t_1-(t_2-t_1)x_1-\ldots-(t_m-t_{m-1})x_m)}\big].
\en
The expression in the right-hand side is invariant under the substitution $t_1\leftrightarrow x_m,\ldots,t_m\leftrightarrow x_1$, which is equivalent to our observation that
the law of ${\cal C}^{(k)}$ is preserved by the reflection about the bisectrix. 

\par 
The $m=1$ instance of (\ref{multdens}) is the 
intensity of ${\cal C}^{(k)}$,
\eq\label{onedens}
p_1(t,x)={k\over \Gamma(k)}(tx)^{k-1}e^{-tx},
\en
which being compared with the (obvious) intensity function $(tx)^{k-1}e^{-tx}/\Gamma(k)$ of the $k$-record process ${\cal R}^{(k)}$
makes us wonder where the factor $k$ is coming from. 
The structure of the processes ${\cal R}^{(k)}$ (which are  neither independent nor identically distributed) is apparently more complex
than that of ${\cal C}^{(k)}$'s. In particular, the process $(R_t^{(k)}, t>0)$ of the value of the last $k$-record observed before $t$ is not 
even Markovian for $k>1$: the law of the life-time at value $x$ is that of the (infinite-mean) random sum $x^{-1}(E_1+\ldots+E_N)$, where all variables are independent,
$E_j$'s are unit exponential, and $N$ is the first success time  in a series of Bernoulli trials with `harmonic' success probabilities  $1/k, 1/(k+1),\ldots$
(explicitly, $\prob(N=n)=1/((n+1)(n+2))$ for $k=2$, but no simple formula exists for $k>2$).   
Still, ${\cal R}^{(k)}$ can be accessed through ${\cal C}^{(k)}$:
\begin{proposition}\label{thin}
The law of ${\cal R}^{(k)}$ is that of a pointwise Bernoulli thinning of ${\cal C}^{(k)}$ with probability $1/k$.
\end{proposition}
\proof Like the thinning argument in Remark \ref{remaf}, this is a consequence of
Lemma \ref{mylemma}.
\endproof
\noindent This explains, of course, the factor $k$ in (\ref{onedens}).

\section{An argument for the ABT lemma}

By Proposition \ref{propo} we can 
identify ${\cal T}^{(k)}$ with $\tau({\cal C}^{(k)})$ and ${\cal S}^{(k)}$ with the set of life-times of $(M_t^{(k)},t>0)$. 
Let ${\cal J}^{(k)}$ be a planar process having points $(s,x)$ where $x\in \xi ({\cal C}^{(k)})$ and $s$ is the life-time
of $(M_t^{(k)},t>0)$ at $x$. It is easily seen that ${\cal J}^{(k)}$ is a marked Poisson process (same as the one denoted 
 $\cal R$ in \cite[p. 47]{ABT}) 
with intensity $(k/x)(xe^{-sx})=ke^{-sx}$; therefore by symmetry of the
intensity   $\tau({\cal J}^{(k)})=_d \xi({\cal J}^{(k)})=_d \xi( {\cal C}^{(k)})$, and  (\ref{ABT}) now follows from
Proposition \ref{propo}.

\par A novelty of this argument is in exploiting the distributional identity of two projections of ${\cal C}^{(k)}$. This allowed us to avoid
 the 
computational part of the proof in \cite{ABT}, where 
one needed to derive the Poisson character of the process ${\cal T}^{(k)}$
from its  
definition  by partial summations over ${\cal J}^{(k)}$.
\par Our proof of (\ref{ABT}) based on the $k$-corners of the unit Poisson process in ${\mathbb R}_+^2$ 
works only for integer $k$. For general $k>0$ one can argue by interpolation from the integer values, since 
all distributions involved depend on the parameter $k$ analytically.
Alternatively, one can introduce $(M_t^{(k)}, t>0)$ for arbitrary $k>0$ directly as a self-similar Markov process entering from the infinity,
as in \cite{Entr}, then derive (\ref{multdens}) and from this conclude about the $t\leftrightarrow x$ symmetry 
of the graph of $t\mapsto M_t^{(k)}$.
By the latter approach one needs to justify (\ref{gamma}), which can be done by application of a moments formula for 
self-similar Markov processes with
the general `stick-breaking' factor like of beta-distributed $B$ above.
\par The case of beta stick-breaking (\ref{beta}) is, in fact, very special 
in that only for this distribution of $B$ the set of jump-times (and the range)
of a self-similar process like $(M^{(k)}_t, t>0)$ is Poisson, see
\cite[Proposition 8]{chain}. The self-similarity property (\ref{selfsim}) persists for arbitrary distribution 
of the stick-breaking factor.

\end{document}